\documentclass[aos,MSNbibl,dvips]{arximspdf}
\usepackage{stfloats}
\usepackage{graphicx}

\doi{10.1214/12-AOS1080} 
\volume{41}
\issue{2}
\pubyear{2013}
\firstpage{436}
\lastpage{463}

\makeatletter

\fnbelowfloat

\newtheorem{prop}{Proposition}[section]
\newtheorem{theorem}[prop]{Theorem}
\newtheorem{cor}[prop]{Corollary}
\newproclaim{defn}[prop]{Definition}
\newproclaim{ex}[prop]{Example}

\newcommand{\independent}{\perp\hspace*{-6.2pt}\perp}

\newcommand{\nVar}{p}
\newcommand{\graph}{G}
\newcommand{\faithPar}{\lambda}
\newcommand{\edgeSet}{E}
\newcommand{\BigSet}{\mathcal{M}_{\graph, \faithPar}}
\newcommand{\vol}{\operatorname{vol}}
\newcommand{\cov}{\operatorname{cov}}
\newcommand{\corr}{\operatorname{corr}}
\newcommand{\PolySet}{\mathcal{P}_{ij\mid S}^{\faithPar}}
\newcommand{\causParJK}{a_{ij}}
\newcommand{\noise}{\varepsilon}
\newcommand{\var}{\operatorname{var}}

\newcommand{\XVect}{X}
\newcommand{\noiseVect}{\varepsilon}

\makeatother

\begin{document}
\begin{frontmatter}

\title{Geometry of the faithfulness assumption in
causal~inference\thanksref{T1}}
\runtitle{Geometry of faithfulness assumption in causal inference}

\thankstext{T1}{Supported in part by NSF Grants
DMS-0907632, DMS-1107000, SES-0835531 (CDI) and ARO Grant
W911NF-11-1-0114; also by the Center
for Science of Information (CSoI), a US NSF Science and Technology
Center, under Grant Agreement CCF-0939370.}

\begin{aug}
\author[A]{\fnms{Caroline} \snm{Uhler}\corref{}\ead[label=e1]{caroline.uhler@ist.ac.at}},
\author[B]{\fnms{Garvesh}~\snm{Raskutti}\ead[label=e2]{raskutti@samsi.info}},
\author[C]{\fnms{Peter}~\snm{B\"uhlmann}\ead[label=e3]{buhlmann@stat.math.ethz.ch}}
\and
\author[D]{\fnms{Bin} \snm{Yu}\ead[label=e4]{binyu@stat.berkeley.edu}}
\runauthor{Uhler, Raskutti, B\"uhlmann and Yu}
\affiliation{IST Austria, SAMSI, ETH Z\"urich and University of
California, Berkeley}
\address[A]{C. Uhler\\
IST Austria\\
3400 Klosterneuburg\\
Austria\\
\printead{e1}}
\address[B]{G. Raskutti\\
SAMSI\\
Research Triangle Park\\
North Carolina 27709\\
USA\\
\printead{e2}\hspace*{47pt}}
\address[C]{P. B\"uhlmann\\
Seminar f\"ur Statistik\\
ETH Z\"urich\\
8092 Z\"urich\\
Switzerland\\
\printead{e3}}
\address[D]{B. Yu\\
Department of Statistics\\
University of California, Berkeley\\
Berkeley, California 94720\\
USA\\
\printead{e4}} 
\end{aug}

\received{\smonth{8} \syear{2012}}
\revised{\smonth{11} \syear{2012}}

%
\begin{abstract}
Many algorithms for inferring causality rely heavily on the
faithfulness assumption. The main justification for imposing this
assumption is that the set of unfaithful distributions has Lebesgue
measure zero, since it can be seen as a collection of hypersurfaces in
a hypercube. However, due to sampling error the faithfulness condition
alone is not sufficient for statistical estimation, and
strong-faithfulness has been proposed and assumed to achieve uniform or
high-dimensional consistency. In contrast to the plain faithfulness
assumption, the set of distributions that is not strong-faithful has
nonzero Lebesgue measure and in fact, can be surprisingly large as we
show in this paper. We study the strong-faithfulness condition from a
geometric and combinatorial point of view and give upper and lower
bounds on the Lebesgue measure of strong-faithful distributions for
various classes of directed acyclic graphs. Our results imply
fundamental limitations for the PC-algorithm and potentially also for
other algorithms based on partial correlation testing in the Gaussian
case.
\end{abstract}

%
\begin{keyword}[class=AMS]
\kwd{62H05}
\kwd{62H20}
\kwd{14Q10}
\end{keyword}
\begin{keyword}
\kwd{Causal inference}
\kwd{PC-algorithm}
\kwd{(strong) faithfulness}
\kwd{conditional independence}
\kwd{directed acyclic graph}
\kwd{structural equation model}
\kwd{real algebraic hypersurface}
\kwd{Crofton's formula}
\kwd{algebraic statistics}
\end{keyword}\vspace*{-2pt}

\end{frontmatter}

\section{Introduction}

Determining causal structure among variables based on observational data
is of great interest in many areas of science. While quantifying
associations among variables is well-developed, inferring causal relations
is a much more challenging task. A popular approach to make the causal inference
problem more tractable is given by directed acyclic graph (DAG) models,
which describe conditional dependence information
and causal structure.\vadjust{\goodbreak}

A DAG $G = (V,E)$ consists of a set of vertices $V$
and a set of directed edges $E$ such that there is no directed cycle. We
index $V=\{1,2,\ldots,p\}$ and consider random variables $\{X_i \mid
i=1,\ldots,p\}$ associated to the nodes $V$. We
denote a directed edge from vertex $i$ to vertex $j$ by $(i,j)$ or
$i\rightarrow j$. In this case, $i$ is called a \textit{parent} of $j$
and $j$ is
called a \textit{child} of $i$. If there is a directed path
$i\rightarrow\cdots\rightarrow j$, then $j$ is called a descendent of
$i$ and $i$ an ancestor of $j$. The skeleton of a DAG $G$ is the
undirected graph obtained from $G$ by substituting directed edges by
undirected edges. Two nodes which are connected by an edge in the
skeleton of $G$ are called \textit{adjacent}, and a triple of nodes
$(i,j,k)$ is an \textit{unshielded triple} if $i$ and $j$ are adjacent
to $k$ but $i$ and $j$ are not adjacent. An unshielded triple $(i,j,k)$
is called a \textit{$v$-structure} if $i\to k$ and $j\to k$. In this case,
$k$ is called a \textit{collider}.

The problem of estimating a DAG from the observational distribution is
ill-posed due to nonidentifiability: in general, several DAGs encode
the same conditional independence (CI) relations and therefore, the true
underlying DAG cannot be identified from the observational
distribution. However, assuming faithfulness (see Definition \ref
{deffaithfulness}), the Markov equivalence class, that is, the
skeleton and the set of $v$-structures of a DAG,
is identifiable (cf.~\cite{pearl00}, Theorem 5.2.6), making it possible
to infer some bounds on causal effects~\cite{makapb09}. We focus here
on the
problem of estimating the Markov equivalence class of a DAG and argue that,
even in the Gaussian case, severe complications arise for data of finite
(or asymptotically increasing) sample size.

There has been a substantial
amount of work on estimating the
Markov equivalence class in the Gaussian case
\cite{spirtes00,chick02,RobScheiSpirtWas,kabu07}. Algorithms which are
based on testing CI relations usually must require
the faithfulness assumption (cf.~\cite{spirtes00}):

\begin{defn}
\label{deffaithfulness}
A distribution $\mathbb{P}$ is \textit{faithful} to a DAG $G$ if no
CI relations other than the ones entailed by
the Markov property are present.
\end{defn}



This means that if a distribution $\mathbb{P}$ is faithful to a DAG
$G$, all
conditional \mbox{(in-)} dependences can be read-off from the DAG $G$ using the
so-called \mbox{$d$-separation} rule (cf.~\cite{spirtes00}). Two nodes
$i,j$ are \textit{$d$-separated} given $S$ if every path between $i$ and
$j$ contains a noncollider that is in $S$ or a collider that is
neither in $S$ nor an ancestor of a node in $S$. For Gaussian models,
the faithfulness assumption can be expressed in terms of the
$d$-separation rule and conditional correlations as follows.
%
\begin{defn}
A multivariate Gaussian distribution $\mathbb{P}$ is said to be
faithful to a DAG $G = (V,E)$ if for any $i,j \in V$ and any $S \subset
V\setminus\{i,j\}$:
\[
j \mbox{ is $d$-separated from } i \mid S\quad\Longleftrightarrow\quad\corr
(X_i,X_j\mid X_S) = 0.
\]
\end{defn}

The main justification for imposing the faithfulness assumption is that
the set of
unfaithful distributions to a graph $G$ has measure zero. However, for
data of finite sample size estimation error issues come into play.
Robins et al.~\cite{RobScheiSpirtWas} showed that many causal
discovery algorithms,
and the PC-algorithm~\cite{spirtes00} in particular, are pointwise
but not uniformly consistent under the faithfulness assumption. This is
because it is
possible to create a sequence of distributions that is faithful but
arbitrarily close to an unfaithful distribution. As a result, Zhang and
Spirtes~\cite{ZhangSpirtes03} defined the strong-faithfulness assumption
for the Gaussian case, which requires sufficiently large nonzero partial
correlations.

%

\begin{defn}
\label{defnfaith}
Given $\lambda\in(0,1)$, a multivariate Gaussian distribution
$\mathbb{P}$ is said to be
\textit{$\lambda$-strong-faithful} to a DAG $G = (V,E)$ if for any
$i,j, \in V$ and any $S \subset
V\setminus\{i,j\}$:
\[
j \mbox{ is $d$-separated from } i \mid S\quad\Longleftrightarrow\quad\bigl|\corr
(X_i,X_j\mid X_S)\bigr| \leq\lambda.
\]
\end{defn}
The assumption of $\lambda$-strong-faithfulness is equivalent to requiring
\[
\min\bigl\{\bigl|\corr(X_i,X_j \mid X_S)\bigr|, j
\mbox{ not $d$-separated from } i \mid S, \forall i,j,S\bigr\} > \lambda.
\]
%
This motivates our next definition which is weaker than
strong-faithfulness.

\begin{defn}
\label{defnrestrfaith}
Given $\lambda\in(0,1)$, a multivariate Gaussian distribution
$\mathbb{P}$
is said to be \textit{restricted $\lambda$-strong-faithful} to a DAG $G =
(V,E)$ if both of the following hold:
\begin{longlist}[(ii)]
\item[(i)] $\min\{| \corr(X_i, X_j\mid X_S)|, (i,j)\in E,
S\subset V\setminus\{i,j\}$ such that $|S|\leq\break\deg(G)\} >
\lambda$,
where here and in the sequel, $\deg(G)$ denotes the maximal degree
(i.e., sum of
indegree and outdegree) of nodes in $G$;
\item[(ii)] $\min\{| \corr(X_i, X_j\mid X_S)|, (i,j,S) \in N_G\} >
\lambda$, where
$N_G$ is the set of triples $(i,j,S)$ such that $i,j$ are not adjacent but
there exists $k\in V$ making $(i,j,k)$ an unshielded triple, and $i$, $j$
are not $d$-separated given $S$.
\end{longlist}
\end{defn}

The first condition (i) is called \textit{adjacency-faithfulness} in
\cite{Zhangtriangle}, the second condition (ii) is called
\textit{orientation-faithfulness}. If a multivariate Gaussian
distribution $\mathbb{P}$ satisfies adjacency-faithfulness with
respect to a DAG $G$, we call the distribution \textit{$\lambda
$-adjacency-faithful} to $G$.
Obviously, restricted
$\lambda$-strong faithfulness is a weaker assumption than
$\lambda$-strong-faithfulness.

We now briefly discuss the relevance of these
conditions and their use in previous work.
Zhang and Spirtes~\cite{ZhangSpirtes03} proved uniform consistency of
the PC-algorithm under the
strong-faithfulness assumption with $\lambda\asymp1/\sqrt{n}$, for the
low-dimensional case where the number of nodes $p = |V|$ is fixed and
sample size
$n \to\infty$. In a high-dimensional and sparse setting, Kalisch and
B\"uhlmann~\cite{kabu07} require strong-faithfulness\vadjust{\goodbreak} with
$\lambda_n \asymp\sqrt{\deg(G)\log(p)/n}$ (the assumption in \cite
{kabu07} is slightly
stronger, but can be relaxed as indicated here). Importantly, since
$\corr(X_i,X_j\mid X_S)$ is required to be
bounded away from $0$ by $\lambda$ for vertices that are not
$d$-separated, the set of distributions
that is not $\lambda$-strong-faithful no longer has measure
$0$.

It is easy to see, for example, from the proof in~\cite{kabu07} that restricted
$\lambda$-strong-faithfulness is a sufficient condition for
consistency of
the PC-algorithm in the high-dimensional scenario [with $\lambda\asymp
\sqrt{\deg(G)\log(p)/n}$] and that the condition is also sufficient and
essentially necessary for consistency of the PC-algorithm. Furthermore,
part (i) of the restricted strong-faithfulness condition is sufficient and
essentially necessary for correctness of the conservative PC-algorithm
\cite{Zhangtriangle}, where correctness refers to the property that
an oriented
edge is correctly oriented but there might be some nonoriented edges which
could be oriented (i.e., the conservative PC-algorithm may not be fully
informative).
The word ``essentially'' above means that we may consider too
many possible separation sets
$S$ where $|S| \le\deg(G)$, while the necessary collection of separating
sets $S$ which the (conservative) PC-algorithm has to consider might be a
little bit smaller. Nevertheless, these differences are minor and we should
think of part (i) of the restricted strong-faithfulness assumption as a
necessary condition for consistency of the conservative PC-algorithm and
both parts (i) and (ii) as a necessary condition for consistency of the
PC-algorithm.

There
are no known upper and lower bounds for the Lebesgue measure of
$\lambda$-strong-unfaithful distributions or of restricted
$\lambda$-strong-unfaithful distributions.
Since these assumptions are so crucial to inferring structure in causal
networks it is vital to understand if restricted and plain
$\lambda$-strong-faithfulness are likely to be satisfied.

In this paper, we address the question of how restrictive the (restricted)
strong-faithfulness assumption is using geometric and combinatorial
arguments. In particular, we develop upper and lower bounds on the
Lebesgue measure of Gaussian distributions that are
not $\lambda$-strong-faithful for various graph structures. By noting
that each
CI relation can be written as a polynomial equation
and the unfaithful distributions correspond to a collection of real
algebraic hypersurfaces, we exploit results from real algebraic
geometry to
bound the measure of the set of strong-unfaithful distributions. As we
demonstrate
in this paper, the strong-faithfulness assumption is restrictive for various
reasons. First, the number of hypersurfaces corresponding to unfaithful
distributions may be quite large depending on the graph
structure, and each hypersurface fills up space in the hypercube. Secondly,
the hypersurfaces
may be defined by polynomials of high degrees depending on the graph
structure. The higher the degree,
the greater the curvature and therefore the surface area of the
corresponding hypersurface. Finally, to get the set of $\lambda
$-strong-unfaithful distributions, these hypersurfaces get fattened up
by a factor which depends on the size of $\lambda$.

Our results show that the set of distributions that do not satisfy
strong-faithfulness can be surprisingly large even for small and sparse
graphs [e.g., 10~nodes and an expected neighborhood (adjacency) size of
2] and small values of $\lambda$ such as $\lambda=0.01$. This implies
fundamental limitations for the PC-algorithm~\cite{spirtes00} and
possibly also for other algorithms based on partial correlations. Other
inference methods, which are not based on conditional independence
testing (or partial
correlation testing), have been described. The penalized maximum
likelihood estimator~\cite{chick02} is an example of such a method and
consistency
results without requiring strong-faithfulness have been given
for the high-dimensional and sparse setting~\cite{geerpb12}. This
method requires, however, a different and so-called permutation
beta-min condition, and it
is nontrivial to understand how the strong-faithfulness condition
and this new condition interact or relate to each other.

The remainder of this paper is organized as follows: Section \ref
{SecMotivatingExample} presents a simple example of a $3$-node fully
connected DAG, where we explicitly list the polynomial equations
defining the hypersurfaces and plot the parameters corresponding to
unfaithful distributions. In Section~\ref{SecProbSetup}, we define the
general model for a DAG on $p$ nodes and give a precise description of
the problem of bounding the measure of distributions that do not
satisfy strong-faithfulness for general DAGs. In Section \ref
{secalgresults}, we provide an algebraic description of the
unfaithful distributions as a collection of hypersurfaces and give a
combinatorial description of the defining polynomials in terms of paths
along the graph. Section~\ref{secbounds} provides a general upper
bound on the measure of $\lambda$-strong-unfaithful distributions and
lower bounds for various classes of DAGs, namely DAGs whose skeletons
are trees, cycles or bipartite graphs $K_{2,p-2}$. Finally, in Section
\ref{secsimulations}, we provide simulation results to validate our
theoretical bounds.

%
\begin{figure}[b]

\includegraphics{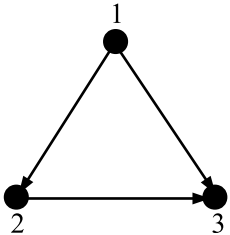}

\caption{Motivating example: 3-node graph.}
\label{fig3nodegraph}
\end{figure}
%

\section{Example: 3-node fully-connected DAG}
\label{SecMotivatingExample}

In this section, we motivate the analysis in this paper using a simple
example involving a $3$-node fully-connected DAG. The graph is shown in
Figure~\ref{fig3nodegraph}. We demonstrate that even in the $3$-node
case, the strong-faithfulness condition may be quite restrictive. We
consider a Gaussian distribution which
satisfies the directed Markov property with respect to the 3-node
fully-connected DAG. An
equivalent model formulation in terms of a Gaussian structural equation
model is given as follows:
\begin{eqnarray*}
X_1 &=& \varepsilon_1,
\\
X_2 &=& a_{12} X_1 + \varepsilon_2,
\\
X_3 &=& a_{13} X_1 + a_{23}
X_2 + \varepsilon_3,
\end{eqnarray*}
where $(\varepsilon_1, \varepsilon_2, \varepsilon_3) \sim
\mathcal{N}(0,I)$.\setcounter{footnote}{1}\footnote{The assumption of $\var(\varepsilon_j)
\equiv1$ is
obviously restricting the class of Gaussian DAG models. We refer to the
more general discussion on this issue in Section \ref
{SecDiscussion}.}\label{foot1}
The parameters $a_{12}, a_{13}$ and
$a_{23}$ reflect the causal structure of the graph. Whether the parameters
are zero or nonzero determines the absence or presence of a directed edge.

It is well known that through observing only covariance information it is
not always possible to infer causal structure. In this example, the
pairwise marginal and the conditional covariances are as follows:%
%
\begin{eqnarray}
\label{eqq1}
\cov(X_1,X_2) & = & a_{12},
\\
\label{eqq2}
\cov(X_1, X_3) & = & a_{13} +
a_{12}a_{23},
\\
\label{eqq3}
\cov(X_2, X_3) & = & a_{12}^2
a_{23} + a_{12}a_{13}+ a_{23},
\\
\label{eqq4}
\cov(X_1, X_2 \mid X_3) & = &
a_{13} a_{23}-a_{12},
\\
\label{eqq5}
\cov(X_1, X_3 \mid X_2) & = & -
a_{13},
\\
\label{eqq6}
\cov(X_2, X_3 \mid X_1) & = & -
a_{23}.
\end{eqnarray}

If it were known a priori that the temporal ordering of the DAG is
$(X_1, X_2, X_3)$, the problem of inferring the DAG-structure would reduce
to a simple estimation problem. We would only need information about the
(non-) zeroes of $\cov(X_1,X_2)$, $\cov(X_1, X_3
\mid  X_2)$ and $\cov(X_2, X_3 \mid  X_1)$, that is,
information whether the
single edge weights $a_{12}, a_{13}$ and $a_{23}$ are zero or not,
which is a standard
hypothesis testing problem. In particular, issues around (strong-)
faithfulness would not arise. However, since
the causal ordering of the DAG is unknown, algorithms based on
conditional independence testing, which
amount to testing partial correlations or conditional
covariances, require that we check \textit{all} partial correlations
between two
nodes given \textit{any subset of remaining nodes}: a prominent example
is the
PC-algorithm~\cite{spirtes00}. For instance for the 3-node case, the
PC-algorithm would
infer that there is an edge between nodes $1$ and $2$ if and only
if $\cov(X_1,X_2) \neq0$ \textit{and} $\cov(X_1, X_2  \mid  X_3)
\neq0$. The issue
of faithfulness comes into play, because it is possible that all causal
parameters $a_{12}, a_{13}$ and $a_{23}$ are nonzero while $\cov(X_1,
X_2 \mid
X_3) = 0$, simply setting $a_{12} = a_{13} a_{23}$ in (\ref{eqq4}).

Since in this example no CI relations are imposed by the Markov property,
a~distribution $\mathbb{P}$ is unfaithful to $G$ if any of the
polynomials in (\ref{eqq1})--(\ref{eqq6}) [corresponding to
(conditional) covariances] are zero. Therefore, the
set of unfaithful distributions for the 3-node example is the union of $6$
real algebraic varieties, namely the three coordinate hyperplanes given
by (\ref{eqq1}), (\ref{eqq5}) and (\ref{eqq6}), two real
algebraic hypersurfaces of degree 2 given by (\ref{eqq2}) and (\ref
{eqq4}), and one real algebraic hypersurface of
degree 3 given by (\ref{eqq3}).

\begin{figure}

\includegraphics{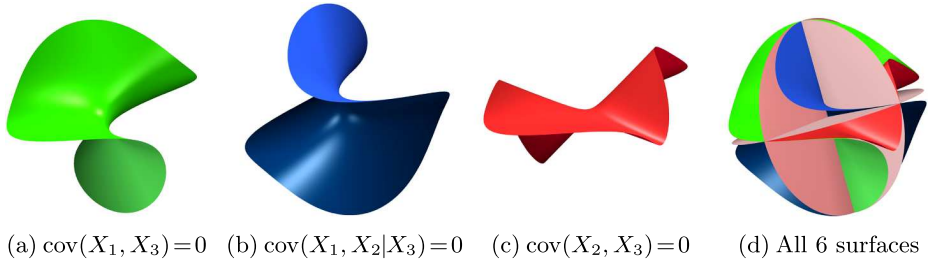}

\caption{Parameter values corresponding to unfaithful distributions in
the 3-node case.}\label{fig3node}
\end{figure}

Assuming that the causal parameters lie in the cube $(a_{12}, a_{13},
a_{23}) \in[-1,1]^3$, we use \texttt{surfex}, a software for visualizing
algebraic surfaces, to generate a plot of the set of parameters leading to
unfaithful distributions. Figure~\ref{fig3node}(a)--(c) shows the
nontrivial hypersurfaces corresponding to $\cov(X_1, X_3) = 0$, $\cov(X_1,
X_2 \mid  X_3) = 0$ and $\cov(X_2, X_3) = 0$. Figure~\ref{fig3node}(d) shows
a plot of the union of all six hypersurfaces.

It is clear that the set of unfaithful distributions has measure
zero. However, due to the curvature of the varieties and the fact that we
are
taking a union of $6$~varieties, the chance of being ``close'' to an
unfaithful distribution is quite large. As discussed earlier, being
close to an unfaithful distribution is of great concern due to sampling
error. Hence, the set of distributions that does not satisfy
$\faithPar$-strong-faithfulness is of interest. As a direct
consequence of Definition~\ref{defnfaith}, this set of distributions
corresponds to the set
of parameters satisfying at least one of the following inequalities:
\begin{eqnarray*}
\bigl|\cov(X_1, X_2)\bigr| & \leq& \faithPar\sqrt{
\var(X_1)\var(X_2)},
\\
\bigl|\cov(X_1, X_3)\bigr| & \leq& \faithPar\sqrt{
\var(X_1)\var(X_3)},
\\
\bigl|\cov(X_2, X_3)\bigr| & \leq& \faithPar\sqrt{
\var(X_2)\var(X_3)},
\\
\bigl|\cov(X_1, X_2 \mid X_3)\bigr| & \leq& \faithPar
\sqrt{\var(X_1\mid X_3)\var (X_2\mid
X_3)},
\\
\bigl|\cov(X_1, X_3 \mid X_2)\bigr| & \leq& \faithPar
\sqrt{\var(X_1\mid X_2)\var (X_3\mid
X_2)},
\\
\bigl|\cov(X_2, X_3 \mid X_1)\bigr| & \leq& \faithPar
\sqrt{\var(X_2\mid X_1)\var(X_3\mid
X_1)}.
\end{eqnarray*}
The set of parameters $(a_{12},a_{13},a_{23})$ satisfying any of the
above relations for $\lambda\in(0,1)$ has nontrivial volume. As we
show in this paper, the volume of the distributions that are not
$\lambda$-strong-faithful grows as the number of nodes and the graph
density grow since both the number of varieties and the curvature of
the varieties increase.

%

\section{General problem setup}
\label{SecProbSetup}

Consider a DAG $G$. Without loss of generality, we assume that the vertices
of $G$ are \textit{topologically ordered}, meaning that $i<j$ for all
$(i,j)\in E$. Each node $i$ in the graph is associated with a random
variable $X_i$. Given a DAG $G$, the random variables $X_i$ are related
to each other by the following structural equations:
%
\begin{equation}
\label{EqnStruct} X_j = \sum_{i < j}{a_{ij}
X_i} + \noise_j,\qquad j = 1,2,\ldots,\nVar,
\end{equation}
where $\noiseVect= (\varepsilon_1,\varepsilon_2,\ldots,\varepsilon_p) \sim
\mathcal{N}(0,I)$
(see footnote~\ref{foot1}) and $a_{ij} \in[-1,+1]$ are the causal parameters
with $a_{ij} \neq0$ if and only if $(i,j) \in\edgeSet$. As we will
see later, we can easily generalize our results to a rescaling of the
parameter cube. In matrix form, these equations can be expressed as
\begin{eqnarray*}
\label{EqnMat} (I - A)^T \XVect= \noiseVect,
\end{eqnarray*}
where $\XVect= (X_1, X_2,\ldots, X_\nVar)$ and $A \in\mathbb
{R}^{\nVar
\times\nVar}$ is an upper triangular matrix with $A_{ij} =
a_{ij}$ for $i < j$. Since $\noiseVect\sim\mathcal{N}(0,I)$,
%
\begin{eqnarray}
\label{EqnDist} \XVect\sim\mathcal{N}\bigl(0, \bigl[(I-A) (I-A)^T
\bigr]^{-1}\bigr).
\end{eqnarray}
We will exploit the distributional form (\ref{EqnDist}) for bounding the
volume of the sets $(a_{ij})_{(i,j) \in\edgeSet} \in
[-1,+1]^{|\edgeSet|}$
that correspond to Gaussian distributions that are not (restricted)
$\lambda$-strong-faithful.

Given $(i,j) \in V\times V$ with $i\neq j$ and $S \subset V \setminus
\{ i,j\}$, we define the set
\begin{eqnarray*}
\label{EqnPoly} \PolySet&:=& \bigl\{ (a_{u,v}) \in[-1,+1]^{|\edgeSet|}
\mid \bigl|\cov (X_i, X_j\mid X_S)\bigr|
\nonumber
\\
&&\hspace*{46pt}\leq \faithPar\sqrt{\var(X_i\mid X_S)
\var(X_j\mid X_S)} \bigr\}.
\end{eqnarray*}
The set of parameters corresponding to distributions that are not
$\lambda$-strong-faithful is
\[
\BigSet:= \mathop{\bigcup_{i,j\in V,  S \subset V\setminus\{
i,j\}:}}_{j\ \mathrm{not}\ d\mbox{-}\mathrm{separated}\ \mathrm{from}\ i \mid  S}
\PolySet.
\]

The set of parameters corresponding to distributions that are not
restricted $\lambda$-strong-faithful
is given by
\[
\label{eqnnecessary} \mathcal{N}_{G,\lambda}^{(1)}:=
\mathop{\bigcup_{i,j\in V,  S
\subset V\setminus\{i,j\}:}}_{(i,j,S)\in N_{G}^{(1)}} \PolySet,
\]
where $N_{G}^{(1)}$ denotes the set of triples $(i,j,S)$,
$S\subset V\setminus\{i,j\}$ with $|S| \le\deg(G)$, satisfying
either\vadjust{\goodbreak} $(i,j) \in E$ or $i$, $j$
are not $d$-separated given $S$ and not adjacent but
there exists $k\in V$ making $(i,j,k)$ an unshielded triple.
The set of parameters corresponding to distributions that are not
$\lambda$-adjacency-faithful [see part (i) of Definition
\ref{defnrestrfaith}] is given by
\[
\label{eqnnecessary2} \mathcal{N}_{G,\lambda}^{(2)}:=
\mathop{\bigcup_{i,j\in V,  S
\subset V\setminus\{i,j\}:}}_{(i,j,S)\in N_{G}^{(2)}} \PolySet,
\]
where $N_{G}^{(2)}$ denotes the set of triples $(i,j,S)$, $S\subset
V\setminus\{i,j\}$ with $|S| \le\deg(G)$, satisfying $(i,j) \in E$.

Our goal is to provide upper and lower bounds on the volume of $\BigSet
$, $\mathcal{N}_{G,\lambda}^{(1)}$ and $\mathcal{N}_{G,\lambda
}^{(2)}$ relative to the volume of $[-1,1]^{|E|}$, that is, to provide
upper and lower bounds for
\[
\frac{\vol(\BigSet)}{2^{|\edgeSet|}}
\quad\mbox{and}\quad
\frac{\vol(\mathcal{N}_{G,\lambda}^{(1)})}{2^{|\edgeSet|}}
\quad\mbox{and}\quad
\frac{\vol(\mathcal{N}_{G,\lambda
}^{(2)})}{2^{|\edgeSet|}}.
\]
This is the probability mass of $\BigSet$, $\mathcal{N}_{G,\lambda
}^{(1)}$ and $\mathcal{N}_{G,\lambda}^{(2)}$ if the parameters
$(a_{ij})_{(i,j) \in
\edgeSet}$ are distributed uniformly in $[-1,+1]^{|\edgeSet|}$, which
we will
assume throughout the paper.



\section{Algebraic description of unfaithful distributions}
\label{secalgresults}


In this section, we first explain that the unfaithful distributions can
always be described by polynomials in the causal parameters
$(\causParJK)_{(i,j)\in E}$ and therefore correspond to a
collection of hypersurfaces in the hypercube $[-1,+1]^{|\edgeSet|}$.
We then
give a combinatorial description of these defining polynomials in terms of
paths in the underlying graph. The proofs can be found in Section~\ref{proofs}.

%
\begin{prop}
\label{proppolynomialCI}
Let $i,j\in V$, $S\subsetneq V\setminus\{i, j\}$ and $Q=S\cup
\{i,j\}$. All CI relations in model (\ref{EqnStruct}) can be
formulated as
polynomial equations in the entries of the concentration matrix
$K=(I-A)(I-A)^T$, namely:


\begin{longlist}
\item $X_i \independent X_j
\Longleftrightarrow (C(K))_{ij}=0$,

\item $X_i \independent X_j \mid X_{V\setminus\{
i,j\}} \Longleftrightarrow K_{ij}=0$,

\item $X_i \independent X_j \mid X_S
\Longleftrightarrow \det
(K_{Q^cQ^c})K_{ij}-K_{iQ^c}C(K_{Q^cQ^c})K_{Q^cj}=0$,
\end{longlist}
where $C(B)$ denotes the cofactor matrix of $B$.\footnote{The
$(i,j)$th cofactor is defined as $C(K)_{ij} = (-1)^{i+j} M_{ij}$ where
$M_{ij}$ is the $(i,j)$th minor of $K$, that is, $M_{ij} = \det(A(-i,-j))$,
where $A(-i,-j)$ is the submatrix of $A$ obtained by removing the $i$th
row and $j$th column of $A$.}
\end{prop}


%


We now give an interpretation of the polynomials defining the hypersurfaces
corresponding to unfaithful distributions in directed Gaussian graphical
models as paths in the skeleton of $G$. The concentration matrix $K$
can be expanded as follows:
\begin{eqnarray*}
K=(I-A) (I-A)^T=I-A-A^T+AA^T.
\end{eqnarray*}
This decomposition shows that the entry
$K_{ij}$, $i\neq j$, corresponds to the sum of all paths from $i$ to
$j$ which lead over a collider $k$ minus the direct path from $i$ to
$j$ if $j$ is a child of $i$, that is,
%
\begin{equation}
\label{eqK} K_{ij}=\sum_{k: i\rightarrow k \leftarrow j}a_{ik}a_{jk}
- a_{ij}.
\end{equation}
Note that $a_{ij}$ is zero in the case that $j$ is not a child of $i$.


For the covariance matrix $\Sigma=K^{-1}$ the equivalent result
describing the path interpretation is given in
\cite{Sullivanttreks}, equation (1), namely
%
\begin{equation}
\label{eqSigma} \Sigma= \sum_{k=0}^{2p-2}
\mathop{\sum_{r+s=k}}_{r,s\leq p-1}
\bigl(A^T\bigr)^rA^s.
\end{equation}
We give a proof using Neumann power series in Section~\ref{proofs}.

Equation (\ref{eqSigma}) shows that the $(i,j)$th entry of $\Sigma$
corresponds to all paths from $i$ to $j$, which first go backwards until
they reach some vertex $k$ and then forwards to $j$. Such paths are called
\textit{treks} in~\cite{Sullivanttreks}. In other words, $\Sigma_{ij}$
corresponds to all collider-free paths from $i$ to $j$.


%
We now understand the covariance between two variables $X_i$ and $X_j$
and the conditional covariance when conditioning on all remaining
variables in terms of paths from $i$ to $j$. In the following, we will
extend these results to conditional covariances between $X_i$ and $X_j$
when conditioning on a subset $S\subsetneq V\setminus\{i, j\}$. This
means that we need to find a path description of
%
\begin{equation}
\label{defpol}
P_{ij\mid S}:= \det(K_{Q^cQ^c})K_{ij}-K_{iQ^c}C(K_{Q^cQ^c})K_{Q^cj}
\end{equation}
[see Proposition~\ref{proppolynomialCI}(iii)] and therefore of the
determinant and the cofactors of $K_{Q^cQ^c}$.

%
%

Ponstein~\cite{Ponstein} gave a beautiful path description of $\det
(\lambda I-M)$ and the cofactors of $\lambda I-M$, where $M$ denotes a
variable adjacency matrix of a not necessarily acyclic directed graph.
By replacing $M$ by $A+A^T-AA^T$, that is by symmetrizing the graph and
reweighting the directed edges, we can apply Ponstein's theorem.



\begin{ponstein*}
Let $i,j\in V$, $S\subsetneq V\setminus\{i, j\}$ and $Q=S\cup\{i,j\}
$ and
let $\hat{G}$ denote the weighted directed graph corresponding to the
adjacency matrix $A+A^T-AA^T$ and $\hat{G}_{Q^c}$ the subgraph resulting
from restricting $\hat{G}$ to the vertices in $Q^c$. Then:
\begin{longlist}[(ii)]
\item[(i)] $\det(K_{Q^cQ^c})=1+\sum_{k=1}^{|Q^c|}  \sum_{m_1+\cdots+m_s=k}(-1)^s \mu(c_{m_1})\cdots\mu(c_{m_s})$,
\item[(ii)] $ (C(K_{Q^cQ^c}) )_{ij}=\sum_{k=2}^{|Q^c|}
\sum_{m_0+\cdots+m_s=k-1}(-1)^s \mu(d_{m_0})\mu(c_{m_1})\cdots\mu
(c_{m_s})$, for $i\neq j$,\vadjust{\goodbreak}
\end{longlist}
where $\mu(d_{m_0})$ denotes the product of the edge weights along a
self-avoiding path from $i$ to $j$ in $\hat{G}_{Q^c}$ of length $m_0$,
$\mu(c_{m_1}), \ldots,\mu(c_{m_s})$ denote the product of the edge
weights along self-avoiding cycles in $\hat{G}_{Q^c}$ of lengths
$m_1,\ldots,m_s$, respectively, and $d_{m_0}, c_{m_1}, \ldots,c_{m_s}$
are disjoint paths.
\end{ponstein*}


Putting together the various pieces in (\ref{defpol}), namely
equation (\ref{eqK}) for describing $K_{QQ}$, $K_{QQ^c}$ and
$K_{Q^cQ}$, and Ponstein's theorem for $\det(K_{Q^cQ^c})$ and
$C(K_{Q^cQ^c})$, we get a path interpretation of all partial correlations.

\begin{ex}
For the special case where the underlying DAG is fully connected and we
condition on all but one variable, that is, $S=V\setminus\{i,j,s\}$, the
representation of the conditional correlation between $X_i$ and $X_j $ when
conditioning on $X_S$ in terms of paths in $G$ is given by
\begin{eqnarray*}
&&
\biggl(1+\sum_{k:  s\rightarrow k}a_{sk}^2
\biggr)
\biggl(\sum_{k:i\rightarrow k\leftarrow j} a_{ik}a_{jk}-a_{ij}
\biggr)
\\
&&\qquad{}- \biggl(\sum_{t:  i\rightarrow t\leftarrow s} a_{it}a_{st}-a_{is}
\biggr) \biggl(\sum_{t:   j\rightarrow
t\leftarrow s}a_{jt}a_{st}-a_{js}
\biggr).
\end{eqnarray*}
\end{ex}

\begin{figure}[b]\vspace*{6pt}

\includegraphics{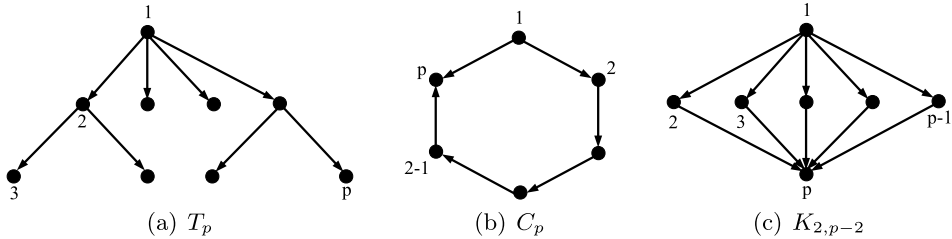}

\caption{Directed tree, cycle and bipartite graph.}
\label{figgraphs}
\end{figure}

In the following, we apply equations (\ref{eqK}), (\ref
{eqSigma}) and Ponstein's theorem to describe the structure of the
polynomials corresponding to unfaithful distributions for various
classes of DAGs, namely DAGs whose skeletons are trees, cycles and
bipartite graphs. We denote by $ T_p$ a directed connected rooted tree
on $p$ nodes, where all edges are directed away from the root as shown
in Figure~\ref{figgraphs}(a). Let $C_p$ denote a DAG whose skeleton is a
cycle, and $K_{2,p-2}$ a DAG whose skeleton is a bipartite graph, where
the edges are directed as shown in Figure~\ref{figgraphs}(b) and
(c).

We denote by $\operatorname{SOS}(a)$ a \textit{sum of squares} polynomial in the
variables\break $(a_{ij})_{(i,j)\in E}$, meaning
\[
\operatorname{SOS}(a)=\sum_k f_k^2(a),
\]
where each $f_k(a)$ is a polynomial in $(a_{ij})_{(i,j)\in E}$. The
polynomials corresponding to unfaithful distributions for the graphs
described in Figure~\ref{figgraphs} are given in the following result.

\begin{cor}
\label{corgraphs}
Let $i,j\in V$ and $S\subset V\setminus\{i,j\}$ such that $i,j$ are
not $d$-separated given $S$. Then the polynomials $P_{ij\mid S}$ defined
in (\ref{defpol}) corresponding to the CI relation $X_i\independent
X_j \mid X_S$ in model (\ref{EqnStruct}) are of the following form:
\begin{longlist}[(c)]
\item[(a)] for $G=T_p$:
\[
a_{i\rightarrow j}\cdot\bigl(1+\operatorname{SOS}(a)\bigr),
\]
where $a_{i\rightarrow j}$ is a monomial and denotes the value of the
unique path from $i$ to $j$;
\item[(b)] for $G=C_p$:
\begin{eqnarray*}
&\displaystyle a_{i\rightarrow j}\cdot\bigl(1+\operatorname{SOS}(a)\bigr)\qquad \mbox{if $ p\notin S$},& \\
&\displaystyle f(\bar{a})a_{i,i+1}-g(\bar{a})a_{j,j+1}\qquad \mbox{if $ S=\{p\}$},&
\end{eqnarray*}
where $a_{i\rightarrow j}$ denotes the value of a path from $i$ to $j$
and $f(\bar{a}), g(\bar{a})$ are polynomials in the variables $\bar
{a}= \{a_{st}\mid(s,t)\notin\{(i,i+1), (j,j+1)\} \}$;
\item[(c)] for $G=K_{2,p-2}$:
\begin{eqnarray*}
&\displaystyle  a_{i\rightarrow j}\cdot\bigl(1+\operatorname{SOS}(a)\bigr)\qquad
\mbox{if $ p\notin S$},& \\
&\displaystyle f(\bar{a})a_{1,j}-g(\bar{a})a_{j,p}\qquad \mbox{if $ i=1$ and $p\in
S$}.&
\end{eqnarray*}
\end{longlist}
\end{cor}

\section{Bounds on the volume of unfaithful distributions}
\label{secbounds}

Based on the path interpretation of the partial covariances explained in
the previous section, we derive upper and lower bounds on the volume of the
parameters that lead to \mbox{$\lambda$-strong-unfaithful} distributions. We also
provide bounds on the proportion of restricted $\lambda$-strong-unfaithful
distributions. These are distributions which do not satisfy the necessary
conditions for uniform or high-dimensional consistency of the
PC-algorithm. Our first result makes use of Crofton's formula for real
algebraic hypersurfaces and the Lojasiewicz inequality to provide a general
upper bound on the measure of strong-unfaithful distributions.

Crofton's formula gives an upper bound on the surface area of a real
algebraic hypersurface defined by a degree $d$ polynomial, namely:

\begin{crofton*}
The volume of a degree $d$ real algebraic hypersurface in the unit
$m$-ball is bounded above by $C(m)d$, where $C(m)$ satisfies
\[
\pmatrix{m+d
\cr
d}-1 \leq C(m) d^{m}.
\]
\end{crofton*}
For more details on Crofton's formula for real algebraic hypersurfaces
see, for example,~\cite{Integergeometry} or~\cite{Guth}, pages 45 and 46.\vadjust{\goodbreak}

The Lojasiewicz inequality gives an upper bound for the distance of a
point to the nearest zero of a given real analytic function. This is
used as an upper bound for the thickness of the fattened hypersurface.
\begin{lojasiewicz*}
Let $f\dvtx  \mathbb{R}^p\rightarrow\mathbb{R}$ be a real-analytic
function and $K\subset\mathbb{R}^p$ compact. Let $V_f\subset\mathbb
{R}^p$ denote the real zero locus of $f$, which is assumed to be
nonempty. Then there exist positive constants $c,k$ such that for all
$x\in K$:
\[
\operatorname{dist}(x,V_f) \leq c\bigl|f(x)\bigr|^k.
\]
\end{lojasiewicz*}

\begin{theorem}[(General upper bound)]
\label{thmupperbound}
Let $G=(V,E)$ be a DAG on $p$ nodes. Then
\begin{eqnarray*}
\frac{\vol(\mathcal{N}_{G,\lambda}^{(2)})}{2^{|\edgeSet|}} &\leq& \frac{\vol(\mathcal{N}_{G,\lambda}^{(1)})}{2^{|\edgeSet|}} \leq \frac{\vol(\BigSet)}{2^{|\edgeSet|}}
\\
&\leq& \frac{C(|E|)c\kappa^k\lambda^k}{2^{{|\edgeSet
|}/{2}}}\sum_{i,j\in V}\sum
_{S \subset V\setminus\{i,j\}}\deg\bigl(\cov (X_i, X_j \mid
X_S)\bigr),
\end{eqnarray*}
where $C(|E|)$ is a positive constant coming from Crofton's formula,
$c, k$ are positive constants, depending on the polynomials
characterizing exact unfaithfulness (for an exact definition, see the
proof), and $\kappa$ denotes the maximal partial variance over all
possible parameter values $(a_{st})\in[-1,1]^{|E|}$, that is,
\[
\kappa=\max_{i,j\in V,  S\subset V\setminus\{i,j\}} \max_{(a_{st})\in[-1,1]^{|E|}} \var(X_i \mid
X_S).
\]
\end{theorem}

Theorem~\ref{thmupperbound} shows that the volume of (restricted)
$\lambda$-strong-unfaithful distributions may be large for two
reasons. First, the number of polynomials grows quickly as the size
and density of the graph increases, and secondly the degree of the
polynomials grows as the number of nodes and density of the graph
increases. The higher the degree, the greater the curvature of the
variety and hence the larger the volume that is filled according to
Crofton's formula. Unfortunately, the upper bound cannot be computed
explicitly, since we do not have bounds on the constants in the
Lojasiewicz inequality.

\begin{pf*}{Proof of Theorem~\ref{thmupperbound}}
It is clear that
\[
\vol\bigl(\mathcal{N}_{G,\lambda}^{(2)}\bigr) \leq \vol\bigl(
\mathcal {N}_{G,\lambda}^{(1)}\bigr) \leq \vol(\BigSet).
\]
Using the standard union bound, we get that
\[
\vol(\BigSet) \leq\mathop{\sum_{i,j\in V,  S \subset
V\setminus\{i,j\}:}}_{j\ \mathrm{not}\ d\mbox{-}\mathrm{separated}\ \mathrm{from}\ i
\mid
S}
{\vol\bigl(\mathcal{P}_{ij\mid S}^{\lambda}\bigr)}.
\]
%
Let $V_{ij\mid S}$ denote the real algebraic hypersurface defined by
$\cov(X_i, X_j \mid X_S)$, that is, the set of all parameter values
$(a_{st}) \in[-1,+1]^{|E|}$ which vanish on $\cov(X_i, X_j \mid X_S)$.
Hence,
\begin{eqnarray*}
\vol\bigl(\mathcal{P}_{ij\mid S}^{\lambda}\bigr) & \leq & \vol\bigl(
\bigl\{ (a_{st}) \in[-1,+1]^{|E|} \mid\bigl|\cov(X_i,
X_j \mid X_S)\bigr|\leq\lambda \kappa\bigr\}\bigr)
\\
& \leq & \vol\bigl(\bigl\{ (a_{st}) \in[-1,+1]^{|E|} \mid
\operatorname{dist} \bigl((a_{st}), V_{ij\mid S} \bigr)\leq
c_{ij\mid S}\lambda^{k_{ij\mid S}} \kappa^{k_{ij\mid S}}\bigr\}\bigr),
\end{eqnarray*}
where $c_{ij\mid S}, k_{ij\mid S}$ are positive constants and the
second inequality follows from the Lojasiewicz inequality.


We apply Crofton's formula on an $|E|$-dimensional ball of radius
$\sqrt{2}$ to get an upper bound on the surface area of a real
algebraic hypersurface in the hypercube $[-1,1]^{|E|}$:
\[
\vol\bigl(\mathcal{P}_{ij\mid S}^{\lambda}\bigr) \leq c_{ij\mid
S}
\lambda^{k_{ij\mid S}} \kappa^{k_{ij\mid S}} 2^{
{|E|}/{2}} C\bigl(|E|\bigr) \deg\bigl(
\cov(X_i, X_j \mid X_S)\bigr).
\]
The claim follows by setting
\[
c=\max_{i,j\in V,  S\subset V\setminus\{i,j\}} c_{ij\mid S}
\quad\mbox{and}\quad k=\min_{i,j\in V,  S\subset
V\setminus\{i,j\}}
k_{ij\mid S}.
\]
\upqed
\end{pf*}

The PC-algorithm in practice only requires $\lambda
$-strong-faithfulness for all subsets $S \subset V\setminus\{i,j\}$
for which $|S|$ is at most the maximal degree of the graph. This could
lead to a tighter upper bound, since we have fewer summands. We will
analyze in Section~\ref{secsimulations} how helpful this is in
practice. In addition, note that we can easily get upper bounds for a
general parameter cube of size $[-r,r]^{|E|}$ by applying Crofton's
formula to a sphere of radius $\sqrt{2}r$.

Since the main goal of this paper is to show how restrictive the
(restricted) strong-faithfulness assumption is, lower bounds on the
proportion of (restricted) $\lambda$-strong-unfaithful distributions are
necessary. However, nontrivial lower bounds for general graphs cannot be
found using tools from real algebraic geometry, since in the worst case the
surface area of a real algebraic hypersurface is zero. This is the case
when the polynomial defining the hypersurface has no real roots. In that
case, the corresponding real algebraic hypersurface is empty. As a
consequence, we need to analyze different classes of graphs separately,
understand the defining polynomials, and find lower bounds for these classes
of graphs. In Section~\ref{secalgresults}, we discussed the
structure of
the defining polynomials for DAGs whose skeleton are trees, cycles or
bipartite graphs, respectively. In the following, we use these results
to find lower
bounds on the proportion of (restricted) $\lambda$-strong-unfaithful
distributions for these classes of graphs.

\begin{theorem}[(Lower bound for trees)]
\label{thmtrees}
Let $T_p$ be a connected directed tree on $p$ nodes with edge set $E$
as shown in Figure~\ref{figgraphs}\textup{(a)}. Then:
\begin{longlist}[(iii)]
\item[(i)]  $\frac{\vol(\mathcal{M}_{T_p,\lambda
})}{2^{|\edgeSet|}}  \geq 1-(1-\lambda)^{p-1}$,\vspace*{2pt}
\item[(ii)]  $\frac{\vol(\mathcal{N}_{T_p,\lambda
}^{(1)})}{2^{|\edgeSet|}}  \geq 1-(1-\lambda)^{p-1}$,\vspace*{2pt}
\item[(iii)]  $\frac{\vol(\mathcal{N}_{T_p,\lambda
}^{(2)})}{2^{|\edgeSet|}}  \geq 1-(1-\lambda)^{p-1}$.
\end{longlist}
\end{theorem}
Theorem~\ref{thmtrees} shows that the measure of restricted and ordinary
$\lambda$-strong-unfaithful distributions converges to
$1$ exponentially in the number $p$ of nodes for fixed
$\lambda\in(0,1)$. Hence, even for trees the strong-faithfulness assumption
is restrictive and the use of the PC-algorithm problematic when the number
of nodes is large.

\begin{pf*}{Proof of Theorem~\ref{thmtrees}}
(i) For a given pair of nodes $i,j\in V$, $i\neq j$, and subset
$S\subset V\setminus\{i,j\}$ we want to lower bound the volume of
parameters $(a_{st})\in[-1,1]^{|E|}$ (in this example $|E|=p-1$) for which
\[
\bigl|\cov(X_i,X_j\mid X_S)\bigr|\leq\lambda\sqrt{
\var(X_i\mid X_S)\var (X_j\mid
X_S)}
\]
or equivalently
\[
|P_{ij\mid S}|\leq\lambda\sqrt{P_{ii\mid S}P_{jj\mid S}}.
\]
From Corollary~\ref{corgraphs}, we know that the defining polynomials
$P_{ij\mid S}$ for $T_p$ are of the form
\[
a_{i\rightarrow j} \cdot\bigl(1+\operatorname{SOS}(a)\bigr).
\]
Similarly as in Corollary~\ref{corgraphs}, one can prove that the
polynomials $P_{ii\mid S}$ are of the form $1+\operatorname{SOS}(a)$ and can therefore
be lower bounded by 1.

So the hypersurfaces representing the unfaithful distributions are the
coordinate planes corresponding to the $p-1$ edges in the tree $T_p$. A
distribution is strong-unfaithful if it is near to any one of the
hypersurfaces (worst case). Since there is a defining polynomial
$P_{ij\mid S}$ without the factor consisting of the sum of squares, the
$\lambda$-strong-unfaithful distributions correspond to the parameter
values $(a_{st})\in[-1,1]^{p-1}$ satisfying
\[
|a_{i\rightarrow j}|\leq\lambda
\]
for at least one pair of $i,j\in V$.
Since we are seeking a lower bound, we set all parameter values to 1
except for one. As a result, a lower bound on the proportion of
$\lambda$-strong-unfaithful distributions is given by the union of all
parameter values $(a_{st})\in[-1,1]^{p-1}$ such that
\[
|a_{st}|\leq\lambda.
\]

We get a lower bound on the volume by an inclusion-exclusion argument.
We first sum over the volume of all by $2\lambda$ thickened coordinate
hyperplanes, subtract all pairwise intersections,\vadjust{\goodbreak} add all three-wise
intersections, and so on. This results in the following lower bound:
\begin{eqnarray*}
\frac{\vol(\mathcal{M}_{T_p,\lambda})}{2^{|\edgeSet|}}&\geq &(p-1)\frac{2\lambda 2^{p-2}}{2^{p-1}}-\pmatrix{p-1
\cr
2}
\frac
{(2\lambda)^2  2^{p-3}}{2^{p-1}}-\cdots
\\
&=&\sum_{k=1}^{p-1}(-1)^{k+1}
\pmatrix{p-1
\cr
k}\lambda^k
\\
&=&1-\sum_{k=0}^{p-1}\pmatrix{p-1
\cr
k}(-
\lambda)^k
\\
&=&1-(1-\lambda)^{p-1}.
\end{eqnarray*}

The proof of (ii) and (iii) is similar. The monomials $a_{i\to j}$
reduce to single parameters $a_{ij}$, since the necessary conditions
only involve $(i,j)\in E$.
\end{pf*}

This theorem is in line with the results in~\cite{Meek}, where they
show that for trees checking if a Gaussian distribution satisfies all
conditional independence relations imposed by the Markov property only
requires testing if the causal parameters corresponding to the edges in
the tree are nonzero.

Note that the behavior stated in Theorem~\ref{thmtrees} is qualitatively
the same as for a linear model $Y=X\beta+\varepsilon$ with active set $S
= \{j\mid \beta_j \neq0\}$. To get consistent estimation of $S$, a
``beta-min'' condition is
required, namely that for some suitable~$\lambda$,
\[
\min_{j \in S}|\beta_j|>\lambda,
\]
meaning that the volume of the problematic set of parameter values
$\beta\in[-1,1]^p$ is given by
\[
1-(1-2\lambda)^{|S|}.
\]
The cardinality $|S|$ is the analogue of the number of edges in a DAG; for
trees, the number of edges is $p-1 \asymp p$ and hence, the comparable
behavior for strong-faithfulness of trees and the volume of coefficients
where the ``beta-min'' condition holds.

Using the lower bound computed in Theorem~\ref{thmtrees}, we can also
analyze some scaling of $n$, $p = p_n$ and $\deg(G) = \deg(G_n)$ as a
function of $n$, such that $\lambda= \lambda_n$-strong-faithfulness
holds. This is discussed in Section~\ref{subsecscalings}.

We now provide a lower bound for DAGs where the skeleton is a cycle on $p$
nodes.
%
\begin{theorem}[(Lower bound for cycles)]
\label{thmcycles}
Let $C_p$ be a directed cycle on $p$ nodes with edge set $E$ as shown
in Figure~\ref{figgraphs}\textup{(b)}. Then:
\begin{longlist}[(iii)]
\item[(i)] $\frac{\vol(\mathcal{M}_{C_p,\lambda
})}{2^{|\edgeSet|}}  \geq 1-(1-\lambda)^{p+{p-1\choose2}}$,\vspace*{2pt}\vadjust{\goodbreak}
\item[(ii)] $\frac{\vol(\mathcal{N}_{C_p,\lambda
}^{(1)})}{2^{|\edgeSet|}}  \geq 1-(1-\lambda)^{3p-2}$,\vspace*{2pt}
\item[(iii)] $\frac{\vol(\mathcal{N}_{C_p,\lambda
}^{(2)})}{2^{|\edgeSet|}}  \geq 1-(1-\lambda)^{2p-1}$.
\end{longlist}
\end{theorem}
For cycles, the measure of $\lambda$-strong-unfaithful distributions
converges to $1$ exponentially in $p^2$. The addition of a single cycle
significantly increases the volume of strong-unfaithful distributions. The
measure of restricted $\lambda$-strong-unfaithful distributions,
however, converges to
$1$ exponentially in $3p$ and hence shows a similar behavior as for
trees. The scaling for achieving strong-faithfulness for cycles is
discussed in
Section~\ref{subsecscalings}.

\begin{pf*}{Proof of Theorem~\ref{thmcycles}}
Similar as for trees, all coordinate hyperplanes correspond to unfaithful
distributions. The corresponding volume of strong-unfaithful distributions
is $2^{p-1}\cdot(2\lambda)$ and there are $p$ such fattened
hyperplanes. In
addition, there are ${p-1\choose2}$ hypersurfaces in the case of (i),
$2(p-1)$ hypersurfaces for (ii), and $p-1$ hypersurfaces for (iii)
defined by polynomials of the form
$f(\bar{a})a_{i,i+1}-g(\bar{a})a_{j,j+1}$, where $\bar{a}= \{
a_{st} \mid(s,t)\notin\{(i,i+1), (j,j+1)\} \}$. Such
hypersurfaces are equivalently defined by
\[
a_{i,i+1}=\frac{g(\bar{a})}{f(\bar{a})}a_{j,j+1}.
\]
Since for any fixed $\bar{a}\in[-1,1]^{p-2}$ this is the
parametrization of a line, we can lower bound the surface area of this
hypersurface by $2^{p-2}\cdot2$, which is the same lower bound as for
a coordinate hyperplane. Similarly as in the proof for trees, an
inclusion-exclusion argument over all hyperplanes yields the proof.
\end{pf*}


Our simulations in Section~\ref{secsimulations} show that by increasing
the number of cycles in the skeleton, the volume of strong-unfaithful
distributions increases significantly. We now provide a lower bound for
DAGs where the skeleton is a bipartite graph $K_{2,p-2}$ and therefore
consists of many 4-cycles. The corresponding scaling for
strong-faithfulness is discussed in Section~\ref{subsecscalings}.
%
\begin{theorem}[(Lower bound for bipartite graphs)]
\label{thmbipartite}
Let $K_{2,p-2}$ be a directed bipartite graph on $p$ nodes with edge
set $E$ as shown in Figure~\ref{figgraphs}\textup{(c)}. Then:
\begin{longlist}[(iii)]
\item[(i)] $\frac{\vol(\mathcal{M}_{K_{2,p-2},\lambda
})}{2^{|\edgeSet|}}  \geq 1-(1-\lambda)^{(p-2)(2^{p-3}+1)}$,\vspace*{2pt}
\item[(ii)] $\frac{\vol(\mathcal{N}_{K_{2,p-2},\lambda
}^{(1)})}{2^{|\edgeSet|}}  \geq
1-(1-\lambda)^{(p-2)(2^{p-3}+1)}$,\vspace*{2pt}
\item[(iii)] $\frac{\vol(\mathcal{N}_{K_{2,p-2},\lambda
}^{(2)})}{2^{|\edgeSet|}}  \geq 1-(1-\lambda)^{(p-2)(2^{p-3}+1)}$.
\end{longlist}
\end{theorem}

\begin{pf}
The graph $K_{2,p-2}$ has $2(p-2)$ edges leading to $2(p-2)$
hyperplanes of surface area $2^{2(p-2)-1}$. In addition, there are
$(p-2)(2^{p-3}-1)$ distinct hypersurfaces defined by polynomials of the
form $f(\bar{a})a_{1,j}-g(\bar{a})a_{j,p}$. Their surface area can be
lower bounded as well by $2^{2(p-2)-1}$ as seen in the proof of
Theorem~\ref{thmcycles}. Hence, the volume of restricted and ordinary
$\lambda$-strong-unfaithful distributions on $K_{2,p-2}$ is bounded
below by
\[
1-(1-\lambda)^{2(p-2)+(p-2)(2^{p-3}-1)}.
\]
\upqed\end{pf}

%
%
%

We remark that we can generalize the lower bounds to a
rescaled parameter cube $[-r,r]^{|E|}$ by replacing $\lambda$ by
$\frac{\lambda}{r}$. Notice that as $r$ increases the lower bounds
decrease but a very large value of $r$ (i.e., very large
absolute values of causal parameters) would be needed to achieve
sufficiently small lower bounds. Furthermore, as discussed in
\cite{singularfaithfulness}, other factors such as singularities on
the partial correlation hypersurfaces may significantly increase the
volume and can occur anywhere on the hypersurface depending on the
structure of the DAG. Therefore, the lower bound may not be tight.

\subsection{Scaling and strong-faithfulness}\label{subsecscalings}

We here consider the setting where the DAG $G = G_n$ and hence the number
of nodes $p = p_n$ and the degree of the DAG $\deg(G) = \deg(G_n)$ depend
on $n$, and we take an asymptotic view point where \mbox{$n \to\infty$}. In such
a setting, we focus on $\lambda= \lambda_n \asymp
\sqrt{\deg(G_n)\log(p_n)/n}$ (see~\cite{kabu07}). We now briefly discuss
when (restricted) $\lambda_n$-strong-faithfulness will asymptotically
hold. For the latter, we must have that the lower bounds (see Theorems
\mbox{\ref{thmtrees}--\ref{thmbipartite}}) on failure of (restricted)
$\lambda_n$-strong-faithfulness tend to
zero.

\textit{Case} I: \textit{lower bound $\asymp1 - (1 - \lambda_n)^{p_n}$.} Such
lower bounds appear for trees (Theorem~\ref{thmtrees}) as well as for
restricted strong-faithfulness for cycles (Theorem~\ref{thmcycles}). The
lower bound $1 - (1 - \lambda_n)^{p_n}$ tends to zero as $n \to\infty
$ if
\[
p_n = o \biggl( \sqrt{\frac{n}{\deg(G_n)\log(n)}} \biggr)\qquad (n \to\infty).
\]
Thus, we have $p_n =
o(\sqrt{n/\log(n)})$ for $\lambda_n$-strong-faithfulness for bounded degree
trees and for restricted $\lambda_n$-strong faithfulness for cycles,
and we have $p_n =
o((n/\log(n))^{1/3})$ for star-shaped graphs.

\textit{Case} II: \textit{lower bound $\asymp1 - (1 - \lambda_n)^{p_n^2}$.}
Such a
lower bound appears for strong-faithfulness for cycles (Theorem
\ref{thmcycles}). The
lower bound $1 - (1 - \lambda_n)^{p_n^2}$ tends to zero as $n \to
\infty$ if
\[
p_n = o \biggl( \biggl(\frac{n}{\deg(G_n)\log(n)} \biggr)^{1/4}
\biggr)\qquad (n \to\infty).
\]
Therefore, we have $p_n =
o((n/\log(n))^{1/4})$ for $\lambda_n$-strong-faithfulness for
cycles.

\textit{Case} III: \textit{lower bound $\asymp1 - (1 - \lambda_n)^{2^{p_n}}$.} This
lower bound appears for strong-faithfulness for bipartite graphs (Theorem
\ref{thmbipartite}). This bound tends to zero as $n \to\infty$ if
\[
p_n = o \bigl(\log(n) \bigr)\qquad (n \to\infty),
\]
regardless of $\deg(G_n)\leq p_n$. Thus, for bipartite graphs with
$\deg(G_n)=p_n-2$ we have $p_n =
o(\log(n))$ for $\lambda_n$-strong-faithfulness.

In summary, even for trees, we cannot have $p_n \gg n$, and
high-dimensional consistency of the PC-algorithm seems
rather unrealistic (unless, e.g., the causal parameters have a distribution
which is very different from uniform).\vspace*{-3pt}

%
%
%

\section{Simulation results}
\label{secsimulations}

In this section, we describe various simulation results to validate the
theoretical bounds described in the previous section. For our simulations,
we used the \texttt{R} library \texttt{pcalg}~\cite{kaetal10}.


In a first set of simulations, we generated random DAGs with a given
expected neighborhood size (i.e., expected degree of each vertex in the
DAG) and
edge weights sampled uniformly in $[-1,1]$. We then analyzed how the
proportion of $\lambda$-strong-unfaithful distributions depends on the
number of nodes $p$ and the expected neighborhood size of the graph.
Depending on the number
of nodes in a graph, we analyzed 5--10 different expected neighborhood
sizes and
generated 10,000 random DAGs for each expected neighborhood
size.

\begin{figure}[b]

\includegraphics{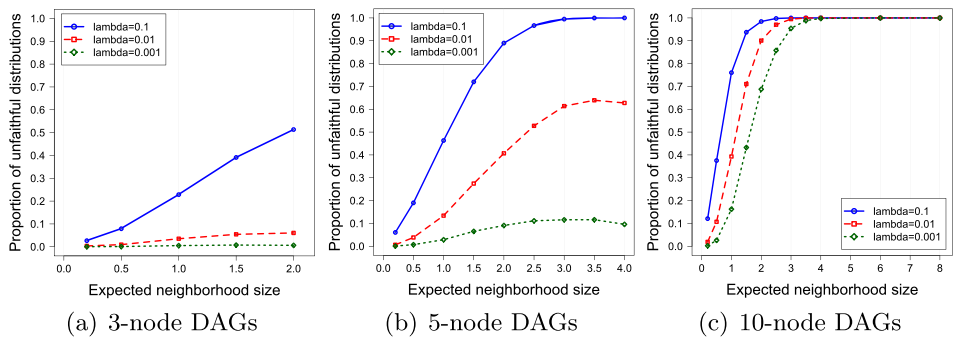}\vspace*{-3pt}

\caption{Proportion of $\lambda$-strong-unfaithful distributions for
3 values of $\lambda$.}
\label{Rplotunfaithful}
\end{figure}

Using \texttt{pcalg} we computed all partial correlations. Since this
computation requires multiple matrix inversions, numerical imprecision has
to be expected. We assumed that all partial correlations smaller than
$10^{-12}$ were actual zeroes and counted the number of simulations, for
which the minimal partial correlation (after excluding the ones with
partial correlation $< 10^{-12}$) was smaller than $\lambda$. The
resulting plots of the proportion of $\lambda$-strong-unfaithful
distributions for three different values of $\lambda$, namely $\lambda=
0.1, 0.01, 0.001$ are given in Figure~\ref{Rplotunfaithful}(a) for $p=3$
nodes, in Figure~\ref{Rplotunfaithful}(b) for $p=5$ nodes and in Figure
\ref{Rplotunfaithful}(c) for $p=10$ nodes.\vadjust{\goodbreak}

It appears that already for very sparse graphs (i.e., expected
neighborhood size of 2) and
relatively small graphs (i.e., 10 nodes) the proportion of
$\lambda$-strong-unfaithful distributions is nearly 1 for $\lambda=0.1$,
about 0.9 for $\lambda=0.01$ and about 0.7 for $\lambda=0.001$. In
addition, the proportion of $\lambda$-strong-unfaithful distributions
increases with graph density and with the number of nodes (even for a
fixed expected neighborhood size). The general upper bound derived in
Theorem~\ref{thmupperbound} shows similar behaviors. The number of
summands and the degrees of the hypersurfaces grow with the number of
nodes and graph density.

\subsection{Bounding the causal parameters away from zero}
\begin{figure}

\includegraphics{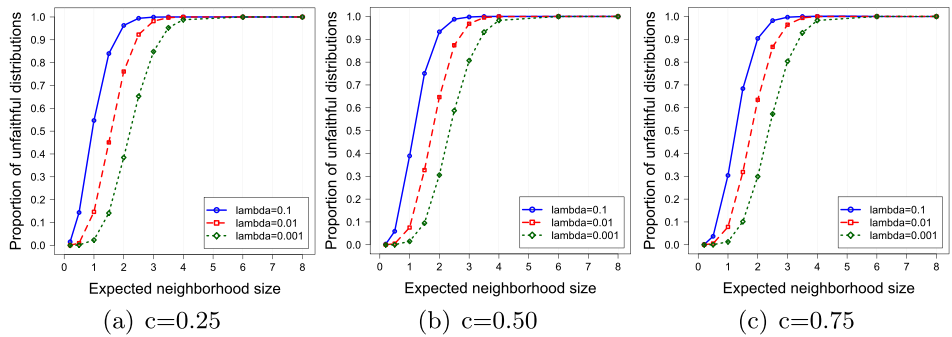}

\caption{Proportion of $\lambda$-strong-unfaithful distributions for
10-node DAGs when restricting the parameter space.}
\label{Rplotrestrictedparameters}
\end{figure}

In the following, we analyze how the proportion of $\lambda
$-strong-unfaithful distributions changes when restricting the
parameter space. The motivation behind this experiment is that
unfaithfulness would not be too serious of an issue if the PC-algorithm
only fails to recover very small causal effects but does well when the
causal parameters are large. We repeated the experiments when
restricting the parameter space to
\[
[-1,-c]\cup[c,1]
\]
for $c=0.25, 0.5$ and 0.75. The results for 10-node DAGs are shown in
Figure~\ref{Rplotrestrictedparameters}. Restricting the parameter
%
space seems to help for sparse graphs but does not seem to play a role
for dense graphs. We now analyze various classes of graphs and their
behavior when restricting the parameter space.

\subsubsection{Trees} 
\begin{figure}

\includegraphics{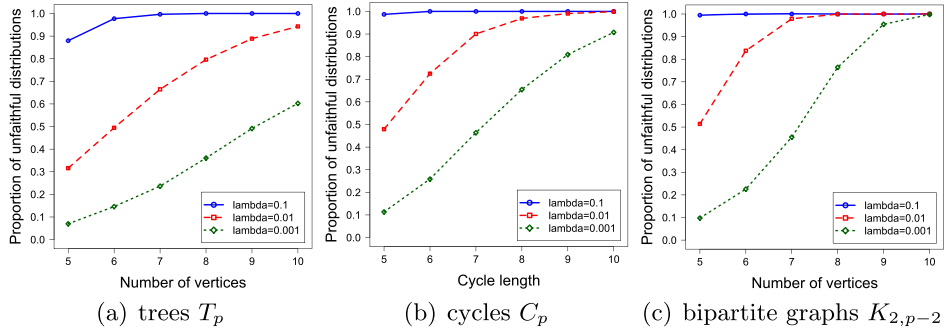}

\caption{Proportion of $\lambda$-strong-unfaithful distributions when
the skeleton is a tree, a cycle or a bipartite graph.}
\label{Rplotdifferentgraphs}
\end{figure}

\begin{figure}[b]

\includegraphics{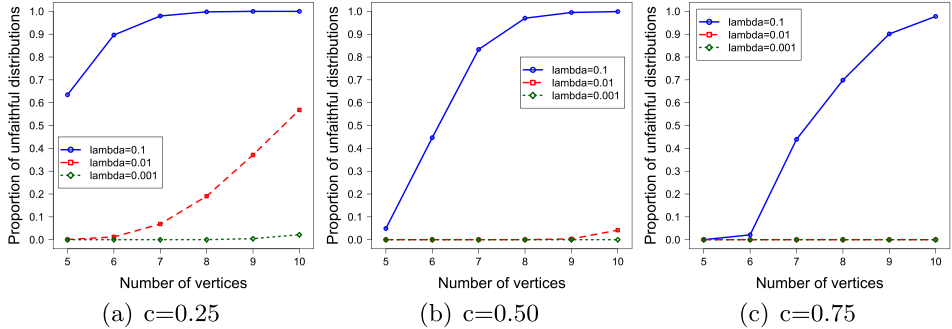}

\caption{Proportion of $\lambda$-strong-unfaithful distributions for
trees when restricting the parameter space.}
\label{Rplottrees}
\end{figure}

We generated connected trees where all edges are
directed away from the root by first sampling the number of levels
uniformly from $\{2,\ldots,p\}$ (a~tree with 2 levels is a star graph, a
tree with $p$ levels is a line), then distributing the $p$ nodes on
these levels such that there is at least one node on each level, and
finally assigning a unique parent to each node uniformly from all nodes
on the previous level. The resulting plots for the whole parameter
space $[-1,1]$ are shown in Figure~\ref{Rplotdifferentgraphs}(a). The plots when
restricting the parameter space for $c=0.25, 0.5$ and 0.75 are shown in
Figure~\ref{Rplottrees}. As before, each proportion is computed from
10,000 simulations.

%
%

For trees restricting the parameter space reduces the proportion of
$\lambda$-strong-unfaithful distributions by a large amount. This can
be explained by the special structure of the defining polynomials
(given in Corollary~\ref{corgraphs}). Since the defining polynomials
of the partial correlation hypersurfaces are of the form
$a_{i\rightarrow j} \cdot(1+\operatorname{SOS}(a))$, the minimal possible value of
these polynomials when restricting the parameter space is
\[
c^{\mathrm{path}\ \mathrm{length}\ \mathrm{from}\ i\ \mathrm{to}\ j}.
\]

\subsubsection{Cycles} 
\begin{figure}

\includegraphics{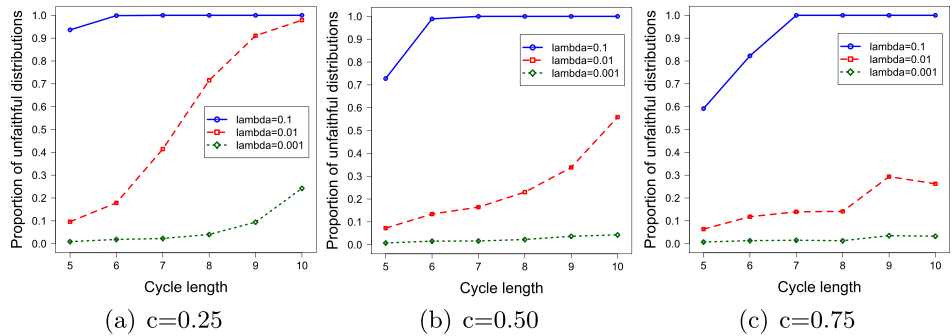}

\caption{Proportion of $\lambda$-strong-unfaithful distributions for
cycles when restricting the parameter space.}
\label{Rplotcycles}
\end{figure}

We generated DAGs where the skeleton is a cycle
and the edges are directed as shown in Figure~\ref{figgraphs}(b). The
edge weights were sampled uniformly\vadjust{\goodbreak} from $[-1,-c]\cup[c,1]$. The
resulting plots for the whole parameter space are shown in Figure
\ref{Rplotdifferentgraphs}(b). The plots for the restricted parameter space with $c=
0.25, 0.5$ and $0.75$ are shown in Figure~\ref{Rplotcycles}. Again,
each point corresponds to 10,000 DAGs.

For cycles restricting the parameter space also reduces the proportion
of $\lambda$-strong-unfaithful distributions, however not as
drastically as for trees. This can again be explained by the special
structure of the defining polynomials (given in Corollary \ref
{corgraphs}). When the defining polynomials are of the form $f(\bar
{a})a_{i,i+1}-g(\bar{a})a_{j,j+1}$, they might evaluate to a very
small number even when the parameters themselves are large.

\subsubsection{Bipartite graphs} 
\begin{figure}

\includegraphics{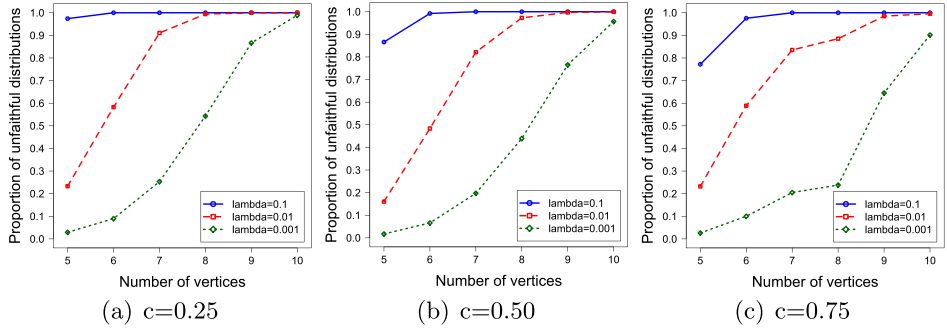}

\caption{Proportion of $\lambda$-strong-unfaithful distributions for
bipartite graphs $K_{2,p-2}$ when restricting the parameter space.}
\label{Rplotbipartite}
\end{figure}

We generated DAGs where the skeleton
is a bipartite graph $K_{2,p-2}$ and the edges are directed as shown in
Figure~\ref{figgraphs}(c). Bipartite graphs $K_{2,p-2}$ consist of
many 4-cycles. For such graphs there are many paths from one vertex to
another and therefore many ways for a polynomial to cancel out, even
when the parameter values are large. As a consequence, for such graphs
restricting the parameter space makes hardly no difference on the
proportion of $\lambda$-strong-unfaithful distributions. This becomes
apparent in Figures~\ref{Rplotdifferentgraphs}(c) and~\ref{Rplotbipartite}.

%
\begin{figure}[b]

\includegraphics{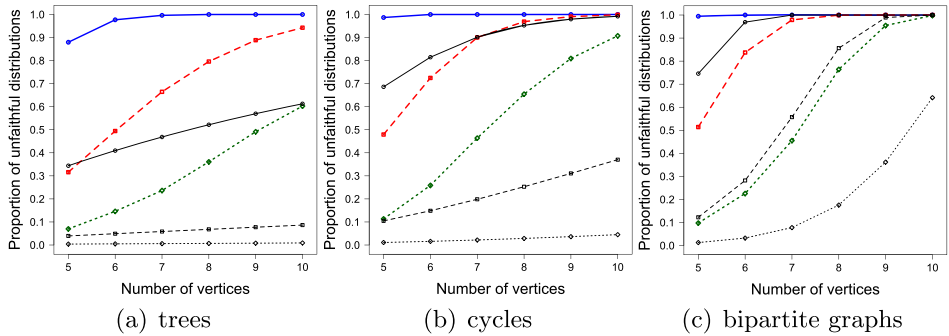}

\caption{Comparison of theoretical lower bounds and approximated
proportion of $\lambda$-strong-unfaithful distributions for trees,
cycles and bipartite graphs $K_{2,p-2}$.}
\label{figbounds}
\end{figure}
%

\subsubsection{Lower bounds} We compare the theoretical lower bounds
derived in Section~\ref{secbounds} to the simulation results in this
section for DAGs where the skeleton is a tree, a cycle or a bipartite
graph when $c=0$. We present our lower bounds together with the
simulation results in Figure~\ref{figbounds}. The black lines
correspond to the lower bounds, the solid line to $\lambda= 0.1$, the
dashed line to $\lambda= 0.01$ and the dotted line to $\lambda
=0.001$. In particular for bipartite graphs our lower bounds
approximate the simulation results very well.

\subsection{\texorpdfstring{Restricted $\lambda$-strong-faithfulness}{Restricted lambda-strong-faithfulness}}
\begin{figure}

\includegraphics{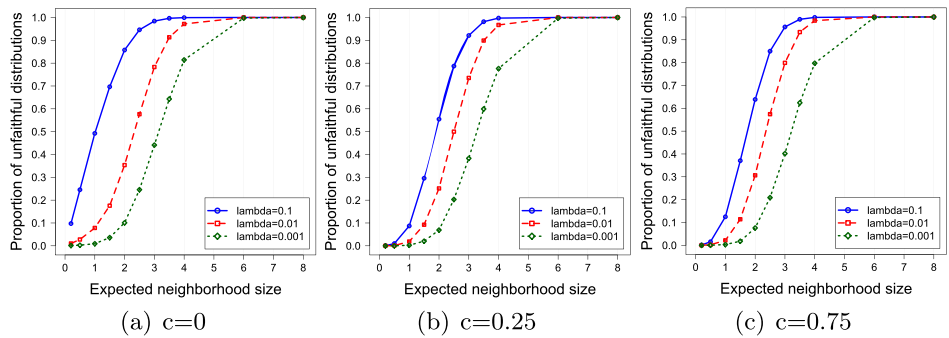}

\caption{Proportion of $\lambda$-adjacency-unfaithful distributions
for 10-node DAGs.}
\label{Rplot10nodesedges}
\end{figure}

As already discussed earlier, the PC-algorithm only requires the
computation of all partial correlations over edges\vadjust{\goodbreak} in the graph $G$ and
conditioning sets $S$ of size at most $\deg(G)$. In order to analyze
when the (conservative) PC-algorithm works, we repeated all our
simulations when restricting the partial correlations to edges in the
graph $G$ and conditioning sets $S$ of size at most $\deg(G)$, that is,
part (i) of the restricted strong-faithfulness assumption in Definition
\ref{defnrestrfaith}, called the adjacency-faithfulness assumption.
The results for general 10-node DAGs are shown in Figure \ref
{Rplot10nodesedges}. We see that the proportion of $\lambda
$-adjacency-unfaithful distributions is slightly reduced compared to
%
the proportion of $\lambda$-strong-unfaithful distributions shown in
Figure~\ref{Rplotrestrictedparameters}, in particular for sparse
graphs. For trees and bipartite graphs the proportion of restricted
$\lambda$-strong-unfaithful distributions is similar to the proportion
of $\lambda$-strong-unfaithful distributions shown in Figures \ref
{Rplotdifferentgraphs},~\ref{Rplottrees} and \ref
{Rplotbipartite}, whereas the behavior for cycles regarding the
proportion of restricted $\lambda$-strong-unfaithful distributions is
similar to trees. We omit these plots here, but remark that they nicely
agree with the theoretical bounds for restricted $\lambda
$-strong-faithfulness and $\lambda$-adjacency-faithfulness derived in
Section~\ref{secbounds}.

%
%
%
%
%
%
%
%
%
%

\section{Discussion}
\label{SecDiscussion}

In this paper, we have shown that the (restricted) strong-faithfulness
assumption is very restrictive, even for relatively small and sparse
graphs. Furthermore, the
proportion of strong-unfaithful distributions grows with the number of
nodes and the number of edges. We have also analyzed the restricted
strong-faithfulness assumption introduced by Spirtes and Zhang \cite
{Zhangtriangle}, a weaker condition than strong-faithfulness, which is
essentially a necessary condition for uniform or high-dimensional
consistency of the popular PC-algorithm and of the conservative
PC-algorithm. As seen in this paper, our lower bounds on restricted
strong-unfaithful distributions are similar to our bounds for strong
faithfulness, implying inconsistent estimation with the PC-algorithm
for a relatively large class of DAGs.

For trees, due to the special structure of the polynomials defining the
hypersurfaces of unfaithful distributions, if the causal parameters are
large, the partial correlations tend to stay away from these
hypersurfaces and
strong-faithfulness holds for a large proportion of distributions.
However, as soon as there are cycles in the graph (even for sparse
graphs), the
polynomials can cancel out also for large causal parameters, and the
strong-faithfulness assumption does not hold. More precisely, if the
skeleton is a single cycle, our lower bounds on the proportion of
restricted strong-unfaithful distributions is of the same order of
magnitude as for trees. However, if the skeleton consists of multiple
cycles as, for example, for bipartite graphs, the lower bounds for
restricted strong-unfaithful distributions are as bad as for plain
strong-unfaithful distributions.

Assuming our framework and in view of the discussion above, in the
presence of
cycles in the skeleton, the (conservative) PC-algorithm is not able to
consistently estimate the true
underlying Markov equivalence class when $p$ is large relative to $n$, even
for large causal parameters (large edge weights). Some special
assumptions on the
sparsity and causal parameters might help, but without making such
assumptions, the limitation is in the range where $p = p_n =
o(\sqrt{n/\log(n)})$. This constitutes a severe
limitation of the
PC-algorithm. As an alternative method, the penalized maximum likelihood
estimator (cf.~\cite{chick02}) does not require strong-faithfulness but
instead a stronger version of a beta-min condition (i.e., sufficiently large
causal parameters)~\cite{geerpb12}. This ``permutation beta-min''
condition has been shown to hold for AR(1) models in~\cite{geerpb12}, page
8. However, a thorough analysis of the ``permutation
beta-min'' condition and a comparison to the strong-faithfulness
condition more generally is quite challenging and remains an
interesting open problem.

Throughout the paper, we have assumed that the causal parameters are
uniformly distributed in the hypercube $[-1,1]^{|E|}$. Since all
hypersurfaces corresponding to unfaithful distributions go through the
origin, a prior distribution which puts more mass around the origin
(e.g., a Gaussian distribution) would lead to a higher proportion of
strong-unfaithful distributions, whereas a prior distribution which
puts more mass on the boundary of the hypercube $[-1,1]$ would reduce
the proportion of strong-unfaithful distributions. Computing and
comparing these measures for different priors would be an interesting
extension of our work. Another interesting problem would be to extend
our results to the case of general error variances [i.e.,
$\operatorname{var}(\varepsilon_j) = \sigma_j^2$]. Finally, very
recently the $k$-triangle-faithfulness assumption has been proposed
\cite{ktrianglefaithfulness} as a sufficient condition for uniform
consistency for inferring certain features of the causal structure.
This assumption is less restrictive than strong-faithfulness, at the
cost of decreasing identifiability, returning a statement
``undecidable'' for some cases. Analyzing how restrictive the
$k$-triangle-faithfulness assumption is and what it means for the
high-dimensional setting represents an interesting future direction.

\section{Proofs}
\label{proofs}

\mbox{}

\begin{pf*}{Proof of Proposition~\ref{proppolynomialCI}}
Statement (i) follows from the matrix inversion formula using the
cofactor matrix, that is,
\[
\Sigma_{ij}=\frac{1}{\det(K)}C(K)_{ij},
\]
and the fact that the concentration matrix $K$ is positive definite and
therefore $\det(K)>0$. Statement (ii) is a well-known fact about the
multivariate Gaussian distribution.

Let $A, B\subset V$ be two subsets of vertices. We denote by $K_{AB}$
the submatrix of $K$ consisting of the entries $K_{ij}$, where
$(i,j)\in A\times B$. Let $K_A$ denote the concentration matrix in the
Gaussian model, where we marginalized over $A^c=V\setminus A$. With
these definitions, we have that
\[
K_A=\Sigma_{AA}^{-1}.
\]

The correlation between $X_i$ and $X_j$ conditioned on $S$ corresponds
to the $(i,j)$th entry in the matrix $K_Q$. Using the Schur complement\vadjust{\goodbreak}
formula, we get that
%
\begin{equation}
\label{eqKSchur} K_Q=K_{QQ}-K_{QQ^c}(K_{Q^cQ^c})^{-1}K_{Q^cQ}.
\end{equation}
Since $K_{Q^cQ^c}$ is positive definite, we can rewrite equation (\ref
{eqKSchur}) as
\[
\det(K_{Q^cQ^c})K_Q=\det(K_{Q^cQ^c})K_{QQ}-K_{QQ^c}C(K_{Q^cQ^c})K_{Q^cQ},
\]
from which statement (iii) follows.
\end{pf*}


\begin{pf*}{Proof of (\ref{eqSigma})}
We first note that the $(i,j)$th element of $A^s$ consists of the sum
of the weights of all paths $p=(p_0,p_1,\ldots, p_s)$ with $p_0=i$ and
$p_s=j$ for which $(p_{k-1},p_{k})\in E$ for all $k=1,\ldots,s$. This
means that $(A^s)_{ij}$ corresponds to all ``forward'' paths from $i$ to
$j$ of length $s$. Analogously, $(A^T)^r$ corresponds to all
``backward'' paths from $i$ to $j$ of length $r$.

We decompose the covariance matrix using the Neumann power series. We
can do this since all eigenvalues of the matrix $A$ are zero (because
$A$ is upper triangular).
\begin{eqnarray*}
\Sigma&=& \bigl((I-A) (I-A)^T \bigr)^{-1}
\\
&=&\sum_{k=0}^{\infty}\sum
_{r+s=k}\bigl(A^T\bigr)^rA^s
\\
&=& \sum_{k=0}^{2p-2}\mathop{\sum
_{r+s=k,}}_{r,s\leq p-1}\bigl(A^T
\bigr)^rA^s.
\end{eqnarray*}
For the last inequality, we used the assumption that the underlying
graph is acyclic. Using the path interpretation it is clear that for
acyclic graphs the matrix $A^s$ is the zero-matrix for all $s\geq p$.
\end{pf*}

\begin{pf*}{Proof of Corollary~\ref{corgraphs}}
To prove (a), we first consider the special case where $G$ is a directed
line on $p$ nodes, where all edges point in the same direction, that
is, $(i,i+1)\in E$ for $1\leq i <p$. The following argument can then
easily be generalized to directed trees $T_p$.

Let $i,j\in V$ and without loss of generality we assume that $i<j$.
Since there are no colliders in $G$, it follows from (\ref{eqK}) that
\[
K_{ij}=\cases{-a_{ij}, &\quad if $j$ is a child of $i$,
\cr
0, &\quad
otherwise,}
\]
$\Sigma_{ij}$ corresponds to all collider-free paths from $i$ to $j$
and therefore
%
\begin{equation}
\Sigma_{ij}= \bigl(1+a_{i-1,i}^2
\bigl(1+a_{i-2,i-1}^2 \bigl(\cdots \bigl(1+a_{12}^2
\bigr) \bigr) \bigr) \bigr)\prod_{k=i}^{j-1}a_{k,k+1}.
\end{equation}
The first term corresponds to the value of all collider-free loops from
$i$ to $i$ and the second term to the value of the path from $i$ to $j$.

Let $S\subsetneq V\setminus\{i, j\}$ and $Q=S\cup\{i,j\}$. If there
exists an element $s\in S$ such that $i<s<j$, then the CI relation $X_i
\independent X_j \mid X_S$ is already entailed by the Markov condition.
We can therefore assume without loss of generality that there is no
$s\in S$ such that $i<s<j$. Since there are no colliders in $G$, it
follows from Proposition~\ref{proppolynomialCI}(iii) that the
corresponding polynomial is of the form
%
\begin{equation}\label{eqpoly1}
\cases{-\det(K_{Q^cQ^c})a_{ij}, &\quad if $j$ is a child of $i$,
\vspace*{2pt}\cr
\displaystyle -\sum_{p,q \in Q^c}a_{ip}C(K_{Q^cQ^c})_{pq}a_{qj},
&\quad otherwise.}
\end{equation}

The corresponding symmetrized and reweighted graph $\hat{G}$ for $p=5$
is shown in Figure~\ref{figchainm}(a). Note that there is a unique
self-avoiding path between any two vertices. As a consequence, the
polynomial corresponding to the CI relation $X_i \independent X_j \mid
X_S$ in (\ref{eqpoly1}) can be written as
%
\begin{equation}
\label{eqpoly2} - \Biggl(1+\sum_{k=1}^{|P|}
\sum_{m_1+\cdots+m_s=k}(-1)^s \mu (c_{m_1})
\cdots\mu(c_{m_s}) \Biggr)\prod_{k=i}^{j-1}a_{k,k+1},
\end{equation}
where $P=Q^c\setminus\{i+1, \ldots,j-1\}$.

\begin{figure}

\includegraphics{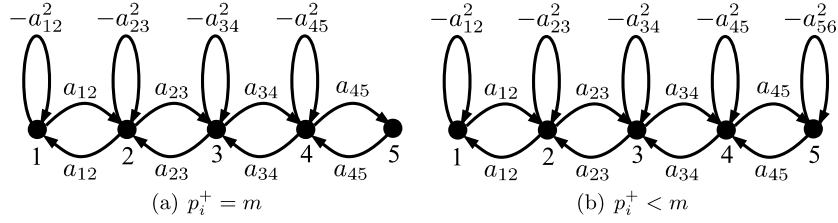}

\caption{Subgraphs $\hat{G}_{P_i}$, where $G$ is a directed line and
$P_i=\{1,2,\ldots,5\}$.}\label{figchainm}
\end{figure}

We now analyze the cycles in $P$. We decompose $P$ into intervals
$P=P_1\cup\cdots\cup P_s$, where $P_i=\{p_i^-, p_i^-+1,\ldots, p_i^+\}
$. We need to distinguish two cases. If $p_i^+=p$, then the subgraph
$\hat{G}_{P_i}$ is of the form as shown in Figure~\ref{figchainm}(a)
(for $p_i^-=1$ and $p_i^+=5$). Otherwise the subgraph is of the form as
shown in Figure~\ref{figchainm}(b) (for $p_i^-=1$ and $p_i^+=5$).

We note that all cycles are either of length 1 (with value
$-a_{k,k+1}^2$) or of length~2 (with value $a_{k,k+1}^2$). In the case
where $p_i^+=p$ all cycles of length 1 cancel with the cycles of length
2. In the case where $p_i^+<p$, however, the cycle of length 1 with
value $-a_{p_i^+, p_i^++1}^2$ does not cancel and therefore neither
does the combination of $k$ cycles
\[
\prod_{j=0}^{k-1}\bigl(-a_{p_i^+-j,p_i^+-j+1}^2
\bigr)
\]
for any $k\in\{1,\ldots,p_i^+-p_i^-\}$. As a consequence, the
polynomial corresponding to the CI relation $X_i \independent X_j \mid
X_S$ in (\ref{eqpoly2}) can be written as
\[
-\prod_{i=1}^s \bigl(1+a_{p_i^+-1,p_i^+}^2
\bigl(1+a_{p_i^+-2,p_i^+-1}^2 \bigl(\cdots \bigl(1+a_{p_i^-,
p_i^-+1}^2
\bigr) \bigr) \bigr) \bigr)\prod_{k=i}^{j-1}a_{k,k+1}.
\]

The proofs for (b) and (c) are analogous and basically require
understanding the cycles in $\hat{G}$.
\end{pf*}

\section*{Acknowledgments}

We wish to thank Marloes Maathuis and Mohab Safey El Din for helpful
discussions. We also thank an Associate Editor and two referees for
constructive comments. Most of the work by the first and second author
was carried out while being at ETH Z\"urich and UC Berkeley,
respectively.



\printaddresses

\end{document}